\documentclass[12pt]{article}
\usepackage{amssymb,amsmath}

\advance \textheight by -1 \topmargin
\advance \textheight by -1 \headheight
\advance \textheight by -1 \headsep
\advance \textheight by -2 \footskip
\vsize = \textheight
\hsize = \textwidth
\def\beql#1#2\eeql{\begin{equation}\label{#1}#2\end{equation}}

\DeclareMathOperator{\Syl}{Syl}

\newtheorem{theorem}{Theorem}[section]

\newtheorem{kor}[theorem]{Corollary}

\newcommand{\bew}{\noindent\underline{Proof.}\ }

\newtheorem{lemma}[theorem]{Lemma}

\newcommand{\be}{\begin{enumerate}}
\newcommand{\ee}{\end{enumerate}}
\newcommand{\bi}{\begin{itemize}}
\newcommand{\ei}{\end{itemize}}

\newcommand{\ba}{\begin{array}}
\newcommand{\ea}{\end{array}}

\newcommand{\disj}{\stackrel{.}{\cup}}

\newcommand{\Z}{{\mathbb{Z}}}
\newcommand{\Q}{{\mathbb{Q}}}
\newcommand{\F}{{\mathbb{F}}}
\newcommand{\N}{{\mathbb{N}}}
\newcommand{\R}{{\mathbb{R}}}

\newcommand{\eb}{\phantom{zzz}\hfill{$\square $}\smallskip}

\renewcommand{\em}{\sf}

\begin{document}

\Large
\begin{center}
{\bf On tight spherical designs.}
\end{center}
\normalsize
\begin{center}
Gabriele Nebe 
\footnote{
Lehrstuhl D f\"ur Mathematik, RWTH Aachen,
Templergraben 64, 52062 Aachen, Germany,
e-mail: nebe@math.rwth-aachen.de
}
 and Boris Venkov
\footnote{
Boris Venkov died in November 2011 before we could finish this paper
}
\end{center}

\small
{\sc Abstract}:
Let $X$ be a tight $t$-design of dimension $n$ for one of the open cases $t=5$ or $t=7$.
An investigation of the lattice generated by $X$ using arithmetic theory of 
quadratic forms allows to exclude infinitely many values for $n$.
\normalsize

\section{Introduction.}

Spherical designs have been introduced in 1977 by 
Delsarte, Goethals and Seidel \cite{DGS} and shortly afterwards 
studied by Eiichi Bannai in a series of papers
(see \cite{Bannai1}, \cite{Bannai2}, \cite{Bannai3} to mention
only a few of them). 
A spherical $t$-design is a finite subset $X$ of the sphere
$$S^{n-1} = \{  x\in \R^n \mid (x,x) = 1 \} $$
such that every polynomial on $\R ^n$ of total degree at most $t$ 
has the same average over $X$ as over the entire sphere. 
Of course the most interesting $t$-designs are those
of minimal cardinality. 
If $t=2m$ is even, then any spherical $t$-design $X\subset S^{n-1}$  satisfies
$$|X| \geq {{n-1+m}\choose{m}} + 
{{n-2+m}\choose{m-1}} $$ and if 
$t=2m+1$ is odd then 
$$|X| \geq 2 {{n-1+m}\choose{m}} .$$ 
A $t$-design $X$ for which equality holds is called 
a {\bf tight} $t$-design.

Tight $t$-designs in $\R^n$ with $n\geq 3$ are very rare.
In \cite{Bannai1} and \cite{Bannai2} it is shown that such tight
designs only exist if $t\leq 5$ and $t=7,11$. 
The tight $t$-designs with $t=1,2,3$ as well as $t=11$ are 
completely classified whereas  their classification for 
$t=4,5,7$ is still an open problem. 
It is known that the existence of a tight $4$-design in dimension 
$n-1$ is equivalent to the existence of a tight $5$-design in
dimension $n$, so the open cases are $t=5$ and $t=7$. 
It is also well known that tight spherical $t$-designs $X$ for odd
values of $t$ are antipodal, i.e. $X=-X$ (see \cite{DGS}).

There are certain numerical conditions on the dimension of such tight
designs. 
A tight $5$-design $X \subset S^{n-1}$ can only exist if 
either $n=3$ and $X$ is the set of 12 vertices of a regular 
icosahedron or
$n=(2m+1)^2-2$ for an integer $m$ (\cite{DGS}, \cite{Bannai1}, \cite{Bannai2}). 
Existence is only known for $m=1,2$ and these designs are unique.
Using lattices \cite{BMV} excludes the next two open cases $m=3,4$ as well as 
an infinity of other values of $m$.
Here we 
exclude infinitely many other cases including $m=6$. 

There are similar results for tight 7-designs. 
Such designs only exist if $n=3d^2-4$ where the only known 
cases are $d=2,3$ and the corresponding designs are unique. 
The paper \cite{BMV} excludes the cases $d=4,5$ and also gives 
partial results on the interesting case $d=6$ which still
remains open.
For odd values of $d$
we use characteristic vectors of the associated odd lattice of odd determinant
to show that $d$ is either $\pm 1 \pmod{16}$ or $\pm 3 \pmod{32}$ 
(see Theorem \ref{oddd}).
We also 
exclude infinitely many even $d$ in 
Theorem \ref{v2deq3}.

\section{General equalities.} 

We always deal with antipodal sets and write them as disjoint union
$$X\disj -X \subset S ^{n-1}(d) = \{ x\in \R^n \mid (x,x) = d\} \mbox{ with } 
s:= |X| \in \N .$$ 
By the theory developed in \cite{Venkov} the set $X\disj -X $ is a 7-design 
 if and only if
 for all $\alpha \in \R^n$
$$ (D6) (\alpha ): \ \ \sum _{x\in X} (x,\alpha )^6 = \frac{3\cdot 5 sd^3}{n(n+2)(n+4)} (\alpha , \alpha )^3.$$
Applying the Laplace operator to $(D6)(\alpha )$  one obtains
$$ (D4) (\alpha ): \ \ \sum _{x\in X} (x,\alpha )^4 = \frac{3sd^2}{n(n+2)} (\alpha , \alpha )^2  \mbox{ and }$$
$$ (D2) (\alpha ): \ \ \sum _{x\in X} (x,\alpha )^2 = \frac{sd}{n} (\alpha , \alpha ) .$$

Substituting $\alpha = \sum _{i=1}^6 \xi _i \alpha _i $ in $(D6)$
$(D4)$ and $(D2)$ we find that for all $\alpha , \beta \in \R^n$:
$$\begin{array}{lll}

(D11) & \sum _{x\in X} (x,\alpha)(x,\beta ) & =  \frac{sd}{n} (\alpha , \beta ) \\
(D13) & \sum _{x\in X} (x,\alpha)(x,\beta )^3 & =  \frac{3sd^2}{n(n+2)} (\alpha , \beta ) (\beta , \beta ) \\
(D22) & \sum _{x\in X} (x,\alpha)^2(x,\beta )^2 & =  \frac{sd^2}{n(n+2)} ( 2 (\alpha , \beta )^2 + (\alpha , \alpha ) (\beta , \beta ))  \\
(D15) & \sum _{x\in X} (x,\alpha)(x,\beta )^5 & =  \frac{3\cdot 5sd^3}{n(n+2)(n+4)}  (\beta , \beta )^2 (\alpha , \beta )  \\
(D24) & \sum _{x\in X} (x,\alpha)^2(x,\beta )^4 & =  \frac{3sd^3}{n(n+2)(n+4)} (  (\beta , \beta )^2 (\alpha,\alpha ) + 4(\alpha , \beta )^2 (\beta , \beta ))  \\
(D33) & \sum _{x\in X} (x,\alpha)^3(x,\beta )^3 & =  \frac{3sd^3}{n(n+2)(n+4)} ( 2 (\alpha , \beta )^3 + 3 (\alpha , \alpha ) (\beta , \beta ) (\alpha , \beta ))  
\end{array}
$$

Similarly $X\disj -X $ is a spherical 5-design, if an only if $(D4)$ and $(D2)$ 
hold for any $\alpha \in \R^n$. Then we obtain the equations $(D11)$, $(D13)$. and $(D22)$.

We will consider the lattice $\Lambda := \langle X \rangle $ 
and $\alpha \in \Lambda ^*$. 
Then $(\alpha , x) $ is integral for all $x\in X$.
This yields certain integrality conditions for the norms and inner
products of elements in $\Lambda ^*$:

\begin{lemma} \label{D4mD2}
If $X\disj -X\subset S^{n-1}(d)$ is a spherical $5$-design then for
all $\alpha, \beta \in \Lambda ^*$
$$  \frac{sd}{12n} (\alpha,\alpha ) (\frac{d}{n+2} (\alpha , \alpha ) -1 ) \in \Z 
$$
and 
$$  \frac{sd}{6n} (\alpha,\beta ) (\frac{d}{n+2} (\alpha , \alpha ) -1 ) \in \Z 
$$
\end{lemma}

\bew
Let $x\in X$ and $k:=(x,\alpha )$.
Then $k^4-k^2$ is a multiple of 12 and hence 
$\frac{1}{12} \sum _{x\in X} (x,\alpha )^4 - (x,\alpha )^2 \in \Z $
which yields the first divisibility condition.
Similarly $k^3-k$ is a multiple of 6 and so
$$\frac{1}{6} \sum _{x\in X} (x,\beta ) ((x,\alpha )^3- (x,\alpha ) ) =
\frac{1}{6} (D13-D11) \in \Z $$
\eb

Similarly 
$$(\beta , x) (\alpha, x) ((\alpha,x) ^2-1)((\alpha , x) ^2 -4) =
(\beta ,x) (\alpha,x)^5-5(\beta ,x) (\alpha ,x)^3+4(\beta,x)(\alpha ,x)$$
is divisible by 5 consecutive integers and hence this quantity is 
a multiple of 120 for any $\alpha,\beta \in \Lambda ^*$ and $x\in X$. 

Moreover $(\alpha, x) ((\alpha,x) ^2-1)$ is divisible by 3 consecutive integers
and therefore a multiple of 6, hence
$$(\beta , x) ((\beta ,x)^2-1) (\alpha, x) ((\alpha,x) ^2-1) =
(\beta ,x) (\alpha , x) ((\beta ,x)^2(\alpha ,x)^2 - (\beta,x)^2-(\alpha,x)^2+1) $$
is divisible by 36. 
Summing over all $x\in X$ we obtain that the right hand side of 
$D15-5D13+4D11$
is a multiple of 120 and that 
$D33-D13-D31+D11 $ is divisible by 36.

\begin{lemma} \label{D15}
If $X\disj -X\subset S^{n-1}(d)$ is a spherical $7$-design then for
all $\alpha, \beta \in \Lambda ^*$
$$\frac{1}{120} 
 (\alpha , \beta ) (\frac{3\cdot 5sd^2}{n(n+2)}(\alpha,\alpha) (
 \frac{d}{n+4}  (\alpha , \alpha ) -1) +
 4 \frac{sd}{n}  )
  \in \Z $$
and
$$\frac{1}{36} (\alpha,\beta) (\frac{3sd^2}{n(n+2)} 
(\frac{d}{n+4} (2(\alpha,\beta)^2+3(\alpha,\alpha)(\beta,\beta)) 
-(\alpha,\alpha) -(\beta,\beta)
) + \frac{sd}{n}  ) \in \Z .$$
\end{lemma}

\section{Tight spherical $7$-designs.}

Let $X\disj -X \subset S^{n-1}(d)$ be a tight spherical $7$-design.
Then $n=3d^2-4$, $(x,y) \in \{ 0,\pm 1 \}$ for all $x\neq y \in X$ and 
$s := |X| = n(n+1)(n+2)/6 .$

Let $\Lambda = \langle X \rangle $ be the lattice generated by 
the set $X$ and put $\Gamma := \Lambda ^*$. 
Then $\Lambda $ is an integral lattice and
$\Lambda $ is even, if $d$ is even. 
Substituting these values into the formulas of Lemma \ref{D15} we obtain

\begin{lemma} \label{div}
For all $\alpha, \beta \in \Gamma $ we have
$$((d^3-d)/240)(\alpha,\beta ) (12d^2 -8 - 15d(\alpha,\alpha) + 5(\alpha,\alpha)^2 ) \in \Z $$
and 
$$((d^3-d)/72) (\alpha,\beta) 
( 3(\alpha,\alpha)(\beta,\beta) - 3d ((\alpha,\alpha) + (\beta,\beta)) 
+ 2 (\alpha,\beta )^2 + (3d^2-2)) \in \Z .$$
\end{lemma}

For a prime $p$ let $v_p$ denote the $p$-adic valuation on $\Q $.

\begin{kor} \label{kdiv} (improvement of \cite[Lemma 4.2]{BMV})
\begin{itemize}
\item[(i)] Let $p \geq 5$ be a prime. If $v_p(d^3-d) \leq 2$ then 
$v_p((\alpha,\alpha))  \geq 0$ for all $\alpha \in \Gamma $. 
\item[(ii)] If $v_3(d^3-d) \leq 4$ then $v_3((\alpha,\alpha )) \geq 0$ 
for all $\alpha \in \Gamma $.
\item[(iii)]
If $v_2(d^3-d) \leq 6$ then $v_2((\alpha ,\alpha )) \geq 0$ for all $\alpha \in \Gamma $. 
\item[(iv)] If $d$ is even but not divisible by $8$ then $v_2((\alpha,\alpha )) \geq 1$ for all
$\alpha \in \Gamma $. 
\item[(v)]
If $d$ is even but not divisible by $32$ then $v_2((\alpha,\beta )) \geq 0 $
for all $\alpha,\beta \in \Gamma $.
\item[(vi)]
If $d$ is odd and $v_2(d^2-1) \leq 4$ then $v_2((\alpha,\beta)) \geq 0 $
for all $\alpha,\beta \in \Gamma $.
\end{itemize}
\end{kor}

\bew
Part (i),(iii) and (iv) are the same as in \cite[Lemma 4.2]{BMV} and follow
from the first congruence in Lemma \ref{div}. 
\\
For (ii) we use the second congruence in the special case $\alpha = \beta $. 
Under the assumption we obtain $v_3((d^3-d)/72) \leq 4-2 \leq 2$. 
If $v_3((\alpha,\alpha )) \leq -1 $ then 
$$v_3 
( 5(\alpha,\alpha)^3 - 6d (\alpha,\alpha)^2 +
 (3d^2-2)(\alpha,\alpha)) = v_3(( \alpha,\alpha)^3) \leq -3$$
contradicting the fact that the product is integral. 
\\
To see (v) we use (iii) to see that $v_2((\alpha,\alpha)) \geq 0$ for 
all $\alpha \in \Gamma $. Then the second congruence yields that 
$v_2(\frac{d}{4}(\alpha,\beta )^3) \geq 0$. Since $v_2(d) < 5$ we obtain $v_2((\alpha,\beta )) \geq 0$.
\\
The last assertion (vi) is obtained by the same argument.
\eb

Using this observation we can extend \cite[Theorem 4.3]{BMV}
which only treats the case
$v_2(d) = 2$.

\begin{theorem} \label{v2deq3}
Assume that 
$v_p(d^3-d) \leq 2$ for all primes $p\geq 5$ and that $v_3(d^3-d)\leq 4$. 
If $v_2(d) = 2,3$ or $4$ then a tight spherical 7-design in dimension
$n=3d^2-4$ does not exist.
\end{theorem}

\bew
$\Gamma $ is integral by 
 Corollary \ref{kdiv} 
 and therefore $\Lambda $ is an even unimodular lattice of
dimension 
$n\equiv 4 \pmod{8}$ which gives a contradiction.
\eb

A similar argument allows to deduce the following lemma from Corollary \ref{kdiv}.

\begin{lemma} \label{oddunimod}
If $d$ is odd and $v_2(d^2-1) \leq 4$ then $\Lambda $ is an 
odd lattice of odd determinant.
If additionally
$v_p(d^3-d) \leq 2$ for all primes $p\geq 5$ and $v_3(d^3-d)\leq 4$ 
then $\Lambda = \Lambda ^*$ is
an odd unimodular lattice.
\end{lemma} 

In particular if $d$ is odd and $d\not\equiv \pm 1 \pmod{16}$ then 
$\Lambda $ is an odd lattice of odd determinant.
Over the 2-adic numbers there is an orthogonal basis 
$$\Lambda \otimes \Z_2 \cong \langle b_1,\ldots , b_n \rangle _{\Z _2} 
\mbox{ with } (b_i,b_j) = 0, (b_k,b_k) = 1, (b_n,b_n) = 1+\delta  \in \{1,3,5,7\}  $$
for $1\leq i \neq j\leq n$, $k=1,\ldots , n-1$. 
Such a lattice contains {\em characteristic vectors}. 
These are elements $\alpha \in \Lambda \otimes \Z_2 $ such that 
$$(\alpha , \lambda ) \equiv (\lambda , \lambda ) \pmod{2} ,\mbox{ for all } 
\lambda \in \Lambda  \otimes \Z_2.$$
Using the basis above, the characteristic vectors in $\Lambda $ are the vectors
$$\alpha = \sum _{i=1}^n a_i b_i \mbox{ with } 
a_i \in 1+2\Z_2 \mbox{ of norm } (\alpha , \alpha ) \equiv n+\delta \pmod{8} .$$ 

\begin{theorem} \label{oddd} 
Let $X\disj -X$ be a tight 7-design 
of dimension $3d^2-4$ with odd $d$.
Assume that $d\not\equiv \pm 1 \pmod{16}$. 
Then either $d\equiv 3 \pmod{32}$ and $\det(\Lambda ) \in (\Z_2^*)^2$ 
or $d\equiv -3 \pmod{32} $ and $\det(\Lambda ) \in 3(\Z_2^*)^2$.
If additionally
$v_p(d^3-d) \leq 2$ for all primes $p\geq 5$ and $v_3(d^3-d)\leq 4$ 
then $d\not\equiv -3 \pmod{16}$.
\end{theorem}

\bew
Let $\Lambda = \langle X \rangle _{\Z_2}$
and $\alpha \in \Lambda $ be a characteristic vector of $\Lambda $ of 
norm $(\alpha,\alpha ) = n+\delta - 8a $ for some $a\in \Z_2$ and $\delta \in \{ 
0,2,4,6 \} $.
Then 
$(\alpha , \lambda ) \equiv (\lambda ,\lambda) \pmod{2} $ for all $\lambda \in \Lambda $, 
in particular $(\alpha,x) $ is odd for all $x\in X$.
For $k>0$ let 
$$n_k := |\{ x\in X \mid (x,\alpha ) = \pm k \} | .$$
Then from (D2), (D4), (D6) we obtain 
$$\begin{array}{lrl} 
(D0) & \sum n_k &  = |X| =  (1/2) (3d^2-4)(3d^2-2)(d^2-1)  \\
(D2) & \sum k^2 n_k &  = (1/2) (3d^2-2)(d^2-1)d (n+\delta -8a ) \\
(D4) & \sum k^4 n_k &  = (3/2) (d^2-1)d^2 (n+\delta -8a)^2  \\
(D6) & \sum k^6 n_k &  = (5/2) (d^2-1) d (n+\delta -8a)^3 . \end{array}$$
Now $n_k\neq 0$ only for odd $k$. 
If $k$ is odd, then $(k^2-1)$ is a multiple of 8 and $(k^2-1)(k^2-9)$ is a multiple 
of $8\cdot 16$. Now $(k^2-1)(k^2-9)(k^2-25) = k^6-35k^4+259k^2-225$ 
is a multiple of $2^{10}3^25$ in particular 
$$(a) \ \ 2^{-7} ((D4)-10(D2)+9(D0)) \in \Z . $$
and 
$$(b) \ \ 2^{-10} ((D6)-35(D4)+259(D2)-225(D0)) \in \Z . $$
We substitute $d=16b+r$ for $r=\pm 3,\pm 5,\pm 7$ into these congruences 
to obtain polynomials in $a$ where the coefficients are polynomials in $b$. 
The contradictions we obtain in the respective cases are listed 
below the table.
$$\begin{array}{|l|c|c|c|c|c|c|}
\hline
r= &                    3 & 5 & 7 & -7 & -5 & -3 \\ 
\hline
\delta = 0 & (c0) &(a2) &(b1) & (a1) & (c2) &(a2) \\
\delta = 2 & (a2) & (c2) & (a1) &(b1) &(a2) &(c0) \\
\delta = 4 & (c1) &(a2) &(b2) & (a1) &(c1) &(a2) \\
\delta = 6 & (a2) &(c1) & (a1) &(b2) &(a2) &(c1) \\
\hline
 \end{array} $$
\begin{itemize}
\item[(a)] In congruence (a) the coefficients of $a$ and $a^2$ are in $\Z[b]$ 
but the constant coefficient is 
\begin{itemize}
\item[(a1)] $p(b) + \frac{b}{2} + \frac{x}{4}$ with 
$p(b) \in \Z[b]$ and $x$ odd.
\item[(a2)] 
 $p(b) + \frac{b}{2} + \frac{x}{8}$ with 
$p(b) \in \Z[b]$ and $x$ odd.
\end{itemize}
\item[(b)]
 In congruence (b) the coefficients of $a$, $a^2$ and $a^3$ are  
in $\Z[b]$ but the constant coefficient is
\begin{itemize}
\item[(b1)] $p(b) +\frac{1}{2} $ with $p(b) \in \Z [b]$.
\item[(b2)] $p(b)+\frac{b}{2}+\frac{x}{4} $ with
$p(b) \in \Z[b]$ and $x$ odd.
\end{itemize} 
\item[(c)] 
 In congruence (b) the coefficient of $a^3$ is in $\Z[b]$ the 
ones of $a$ and $a^2$ 
are in $\frac{1}{2} + \Z[b]$ but the constant coefficient is
\begin{itemize}
\item[(c0)] $p(b) + \frac{b}{2}$ with $p(b) \in \Z[b] $. 
Here we can only deduce that $b$ is even.
\item[(c1)] $p(b) + \frac{x}{8}$ with $p(b)\in \Z[b]$ and $x$ odd. 
\item[(c2)] $p(b) + \frac{b}{2} + \frac{x}{4}$ with $p(b)\in \Z[b]$ and $x$ odd.
\end{itemize}
\end{itemize}
Hence only the cases $r=3$, $\delta = 0$ and  $r=-3$, $\delta =2$ 
are possible and then  $b$ is even. 
\eb

To summarize we list a few small values that are excluded by Theorem \ref{oddd} and Theorem \ref{v2deq3}:
\begin{kor}
There is no tight 7-design of dimension $n=3d^2-4$ for 
$$d\in \{ 4,5,7,8,9,11,12,13,16,19,20,21,\ldots   \} $$
\end{kor}

\section{Tight spherical $5$-designs.}

Assume that $d=2m+1$ and that $X\disj -X $ is a tight spherical 
5-design in dimension $n=d^2-2$.
Then $|X| = n(n+1)/2$ and scaled such that 
$(x,x) = d$ for all $x\in X$ we have $(x,y) = \pm 1$ for $x\neq y \in X$ 
and $\Lambda := \langle X \rangle $ is an odd integral lattice.
With these values the formula $(D4)$ reads as
$$
(D4) \ \  \sum _{x\in X} (x,\alpha )^4 = 6m(m+1) (\alpha , \alpha )^2  .
$$

\begin{lemma} (see \cite[Lemma 3.6]{BMV}) 
Assume that  $m(m+1)$ is not divisible by the square of a 
prime $p \geq 5$. Then $(\alpha, \alpha ) \in \Z[1/6] $ for 
all $\alpha \in \Lambda ^*$.
\end{lemma} 

Substituting the special values into the 
formula of  Lemma \ref{D4mD2}  we immediately obtain

\begin{lemma} (see \cite[Lemma 3.3]{BMV}) \label{div3}
For all $\alpha \in \Lambda ^*$
$$\frac{1}{6} m(m+1)(\alpha ,\alpha ) ( 3(\alpha, \alpha ) - (2m+1)) \in \Z $$
\end{lemma}

\begin{kor} \label{even*}
If $m(m+1) $ is not a multiple of $8$, then 
$(\alpha ,\alpha ) \in \Z_2$ is $2$-integral for all $\alpha \in \Lambda ^*$.
\end{kor}

We now treat the Sylow 3-subgroup $D_3:= \Syl_3(\Lambda^*/\Lambda )$.

\begin{lemma} \label{D3}
Assume that $m(m+1)$ is not a multiple of $9$.
Then $|D_3| \in \{ 1, 3 \}$.
\end{lemma}

\bew
Assume that $D_3\neq 1$.
Since $D_3$ is a regular quadratic $3$-group it contains an 
anisotropic element $\alpha +\Lambda \in \Lambda ^*/\Lambda $
with $(\alpha,\alpha ) = \frac{p}{q}$ and $3\mid q$.
By equality $(D4)$ the denominator $q$ is not divisible by $9$, 
in particular the exponent of $D_3$ is 3 and 
 $(\alpha , \alpha ) = \frac{p}{3}$ with a 3-adic unit $p
\equiv \pm 1 \pmod{3}$.
Now Lemma \ref{div3} gives
$$\frac{1}{18} m(m+1)p ( p - (2m+1)) \in \Z $$
Since $m(m+1)$ is not a multiple of 9, this implies that 
$p\equiv (2m+1) \pmod{3}$. 
If $|D_3| > 3$, then the regular quadratic $\F_3$-space 
$D_3$ is universal, representing also elements 
$\frac{p}{3}$ with $p\not\equiv (2m+1) \pmod{3}$.
This is a contradiction. 
So $|D_3| = 1$ or $|D_3| = 3$.
\eb

Let $\Lambda _+ $ be the even sublattice of
$\Lambda = \langle X \rangle $. 
Then $\Lambda = \Lambda _+ \disj \Lambda _-$ with
$\Lambda _- = x + \Lambda _+$ for any $x\in X$.
Since $(x,y)$ is odd for all $x\in X$, the lattice
$$\Lambda _+ = \{ \sum _{x\in X} c_x x \mid c_x \in \Z, \sum _{x\in X} c_x 
\mbox{ even } \} $$ and
$(\alpha , x) \in 2\Z $ for any $\alpha \in \Lambda _+$ and $x\in X$.
Therefore $\Lambda _+ \subset 2 \Lambda ^*$ and the 
lattice $\Gamma := \frac{1}{\sqrt{2}} \Lambda _+$ is an integral 
lattice of dimension $n$. 

The next lemma is an improvement of \cite[Lemma 3.6]{BMV}.

\begin{lemma} \label{C2}
Assume that $m(m+1)$ is not divisible by the square of an odd prime and 
that $m$ is odd and $(m+1)$ is not a multiple  of $8$.
Then for any $x\in X$
$$\Gamma^*/\Gamma = \langle \frac{1}{\sqrt{2}} x + \Gamma \rangle
\cong \Z /2\Z .$$
\end{lemma}

\bew
For odd primes $p$ the Sylow p-subgroup of $\Gamma ^*/\Gamma $ is
isomorphic to the one of $\Lambda ^*/\Lambda $ and hence 
$\{0 \}$ for $p\geq 5$ and either $\{0 \}$ or $\Z/3\Z $ for $p=3$.
Clearly $\alpha := \frac{1}{\sqrt{2}} x  \in \Gamma ^*$ has order 2
modulo $\Gamma $.
Moreover 
$$\Gamma ^*  = \sqrt{2} \Lambda _+^* =  \langle \alpha , \sqrt{2} \Lambda ^*
\rangle $$ 
is an overlattice of $\sqrt{2}\Lambda ^*$ of index 2. 
Now by Corollary \ref{even*} 
$(\beta , \beta ) \in 2\Z_2$ for all elements $\beta \in \sqrt{2}\Lambda ^*$
and since $x\in \Lambda $ we get $(\beta,\alpha) \in \Z $ for all
$\beta \in \sqrt{2} \Lambda ^*$. 
Since the Sylow 2-subgroup $D_2$ of $\Gamma ^*/\Gamma $ is a regular
quadratic $2$-group and $D_2 \cap \sqrt{2}\Lambda ^*/\Gamma $ is in 
the radical of this group we obtain that $D_2=\langle \alpha + \Gamma \rangle
\cong \Z/2\Z $.
To exclude the case that $D_3 = \Z/3\Z $ we use the fact that 
$\Gamma $ is an even lattice and hence the 
Gau\ss \ sum 
$$G(\Gamma ):= 
\frac{1}{\sqrt{2\cdot 3^t}} \sum _{d\in \Gamma^*/\Gamma } \exp (2\pi i q(d))$$
for the quadratic group 
$(\Gamma ^*/\Gamma, q)$ with $q(z+\Gamma ):= \frac{1}{2}(z,z) + \Z $
equals
$$G(\Gamma ) = \exp(\frac{2\pi i}{8}) ^{n} 
=\exp(\frac{2\pi i}{8}) ^{-1}  $$
by the Milgram-Braun formula.
Clearly $G(\Gamma )$
is the product of the Gau\ss \ sums of its Sylow subgroups,
$G(\Gamma ) = G_2 G_3$ with 
$$G_2 = \frac{1}{\sqrt{2}} (1+\exp(2\pi i \frac{2m+1}{4})) = \frac{1-i}{\sqrt{2}}
=\exp(\frac{2\pi i}{8}) ^{-1} = G(\Gamma ) $$
since $m$ is odd.
This implies that $G_3 = 1$.
Then \cite[Corollary 5.8.3]{Scharlau} shows that 
$D_3$ cannot be anisotropic, and hence by Lemma \ref{D3} 
$D_3=\{ 0 \}$.
\eb

\begin{theorem} (see also \cite[Theorem 3.10]{BMV} for one case)
Assume that $m(m+1)$ is not divisible by the square of an odd prime,
$m$ is even but not divisible by $8$.
Then $\Gamma ^*/\Gamma \cong \Z/6\Z $ and $m\equiv -1 \pmod{3}$.
\end{theorem}

\bew
With the same proof as above we obtain $G(\Gamma )= \exp(\frac{2\pi i}{8}) ^{-1}$
and $G_2 = \exp(\frac{2\pi i}{8}) $ and
hence $G_3=-i$.
Then \cite[Corollary 5.8.3]{Scharlau} yields that 
$D_3 = \langle \beta + \Gamma \rangle $ with 
$3(\beta ,\beta ) \equiv 1 \pmod{3} $.
Let $\lambda := \sqrt{2} \beta \in \Lambda ^*$.
Then $(\lambda ,\lambda ) = \frac{p}{3}$ with $p\equiv 2 \pmod 3$. 
Then the integrality condition in Lemma \ref{div3} shows that 
$$ m(m+1)(2m-1) \in 9 \Z_3 $$
is a multiple of $9$. 
This implies that 
$m\not\equiv 1 \pmod{3} $ as it was already observed in \cite{BMV}
but also that $m\not\equiv 0 \pmod{3}$. 
\eb

\begin{kor}
$m\neq 3,4,6,10,12,22,28,30,34,42,46,\ldots $.
\end{kor}


%
%

\end{document}